\documentclass[12pt]{article} 
\newdimen\paperhight
\usepackage{amsmath}
\usepackage{amssymb}
\usepackage{amsfonts}
\usepackage{graphicx}
\usepackage[usenames]{color}
\usepackage{hyperref}
\hypersetup{
  colorlinks   = true, %Colours links instead of ugly boxes
  urlcolor     = blue, %Colour for external hyperlinks
  linkcolor    = blue, %Colour of internal links
  citecolor   = red %Colour of citations
}

\numberwithin{equation}{section} %Changes how the equation labeling works (commenting out turns equation labels to single numbers.
\newcommand{\SL}{{\rm SL}}
\newcommand{\pd}{\partial} 

\newcommand{\pr}{\par \vspace{3mm}\noindent [{\bf Proof}] \qquad}
\newcommand{\prend}{\hfill \qed \par \vspace{3mm}}
\newcommand{\qed}{\quad\hbox{\rule[-2pt]{3pt}{6pt}}\par\vspace{3mm}}

\newcommand{\1}{{\bf 1}} 
\newcommand{\C}{\mathbb C} 
\newcommand{\Z}{\mathbb Z} 
\newcommand{\Q}{\mathbb Q} 
 
\newcommand{\N}{\mathbb N} 
\newcommand{\R}{\mathbb R}

\newcommand{\CG}{{\cal G}}

\newcommand{\CH}{{\cal H}}

\newcommand{\tr}{{\rm tr}}
\newcommand{\wt}{{\rm wt}}
\newcommand{\Tr}{{\rm Tr}}
\newtheorem{thm}{Theorem}
\newtheorem{prn}[thm]{Proposition}

\newtheorem{lmm}[thm]{Lemma}
\newtheorem{rmk}[thm]{Remark}

\begin{document}
\title{A modular invariance property of multivariable trace functions 
for regular vertex operator algebras}

\author{\begin{tabular}{c}
Matthew Krauel\footnote{e-mail: matthew.krauel@gmail.com. Supported by the Japan Society of the Promotion of Science (JSPS), No. P13013,
and the European Research Council (ERC) Grant agreement n. 335220 - AQSER.} \\
{\small Mathematical Institute, University of Cologne, Germany} \\ 
{\small and} \\
Masahiko Miyamoto\footnote{e-mail: miyamoto@math.tsukuba.ac.jp. Partially supported 
by the Grants-in-Aids for Scientific Research, 
No. 22654002, The Ministry of Education, Science and 
Culture, Japan.} \\
 {\small Institute of Mathematics, University of Tsukuba, Japan}
\end{tabular}}
\date{}
\maketitle

\begin{abstract}
We prove an $\SL_2 (\Z)$-invariance property of multivariable trace functions 
on modules for a regular VOA. 
Applying this result, we provide a proof of the inversion transformation formula for 
Siegel theta series. As another application, we show that if $V$ is a simple regular VOA containing a simple regular subVOA $U$ whose commutant $U^c$ is simple, regular, and satisfies $(U^c)^c =U$, then all simple $U$-modules appear in some simple $V$-module.
\end{abstract}

\renewcommand{\baselinestretch}{1.5}

\section{Introduction} \label{Section-Intro}

The concept of a vertex operator algebra (VOA) was introduced by Borcherds \cite{B} to explain 
a mysterious relation between the Monster simple group and the elliptic modular function $J(\tau)$.
In the years since, this connection has been elucidated further and generalized to encompass a wide class of VOAs and elliptic modular forms. In the heart of this developing theory reside trace functions over modules of endomorphisms associated with the VOA. In particular, these functions include an operator formed from a matching of a distinguished element from the VOA, and a single variable in the complex upper half-plane. Meanwhile, the element resulting from this pairing resides in a one-dimensional Jordan subalgebra of the VOA, and begs the question whether trace functions exist which instead incorporate elements from larger Jordan subalgebras. The primary aim of this paper is to study such multivariable trace functions and establish functional equations for them with respect to the group $\text{SL}_2 (\mathbb{Z})$.

The development of these equations utilizes a seminal result of Zhu \cite{Z}, which establishes that the space of trace functions on simple modules of a regular VOA is invariant under the standard action of $\text{SL}_2 (\mathbb{Z})$. In particular, Zhu shows that the action of an element of $\text{SL}_2 (\mathbb{Z})$ on a single-variable trace function on a simple module is a linear combination of the trace functions for all simple modules of the VOA with coefficients dependent on the representation of the element in $\text{SL}_2(\mathbb{Z})$. As we lift Zhu's theory to the multivariable case below, we find that we recover these same coefficients. Using Verlinde's formula, we exploit this fact to show that every simple module of a regular subVOA whose commutant satisfies certain conditions is contained in a simple module of the VOA (see Theorem \ref{MC} below).

Beyond considering such regular subVOAs and their commutants, a number of important classes of VOAs are known to contain appropriate Jordan subalgebras and fit the framework presented here to construct multivariable trace functions. We discuss some of these below and look more closely at an application to lattice VOAs, where we formulate another proof of the transformation properties for Siegel theta functions. To explain our results in more detail, we first review the relevant theory and notation pertaining to VOAs.

A VOA is a quadruple $(V,Y(\cdot,z),\1,\omega)$, which we simply denote by $V$,  
consisting of a graded vector space $V=\oplus_{n\in \Z} V_n$, 
a linear map $Y(\cdot,z):V \to {\rm End}(V)[[z^{-1},z]]$, and 
two notable elements $\1\in V_0$ and 
$\omega\in V_2$ called the Vacuum and Virasoro elements, respectively. We say $v$ has weight $n$ if $v\in V_n$ and denote the weight of $v$ by $\wt(v)$ if it is not specified.
An image $Y(v,z)=\sum_{n\in \Z} v_nz^{-n-1}$ of $v\in V$ is called a vertex operator of $v$, and it can be shown that $v_{\wt(v)-1}$ is a weight-preserving operator for a homogeneous element 
$v$. We denote this unique operator by $o(v)$ and extend it linearly. Meanwhile, the operators
$L(n)$ defined by $Y(\omega,z)=\sum_{n\in \Z}L(n)z^{-n-2}$ satisfy a Virasoro algebra bracket relation 
$$ [L(n),L(m)]=(n-m)L(n+m)+\delta_{n+m,0}\frac{n^3-n}{12}c$$
for some $c\in \C$, called the central charge of $V$. The eigenvalues of $L(0)$ provide the weights on $V$, that is, $V_n =\{ v\in V \mid L(0)v=nv \}$. 

In this paper, we assume that $V$ is a regular VOA of CFT-type 
(i.e. $C_2$-cofinite, rational, and $\N$-graded with $V_0=\C \1$) of central charge $c$. A number of important consequences can be drawn from these assumptions. For one, such a $V$ has only finitely many isomorphism classes of simple $V$-modules $\{W^1,\ldots,W^r\}$, and all of them are $\N$-gradable. We refer the reader to \cite{DLM-Reg} for a further discussion on the definition of regular VOAs and implications of this definition on the structure of $V$-modules. Another consequence of a regular VOA which stems from the CFT-type assumption, is a refinement on the classification of symmetric invariant bilinear forms on $V$. Indeed, in \cite{Li} it is shown that the space of bilinear forms for such a VOA is isomorphic to the dual of $\left(V_0/L(1)V_1\right)$. Therefore, if $L(1)V_1=0$, the choice of an $\alpha \in \mathbb{C}$ satisfying $\langle \textbf{1},\textbf{1}\rangle =\alpha$ uniquely defines a bilinear form on $V$. If $\alpha$ is chosen to be $-1$, then for any $a,b\in V_1$ we have $a_1 b=\langle a,b\rangle \textbf{1}$. See \cite{Li} for details and proofs concerning bilinear forms on $V$.

Additionally, for a regular VOA $V$ we may invoke Zhu's \cite{Z} (see also \cite{DLM-Orbifold}) modular-invariance results for single-variable trace functions mentioned above. Specifically, Zhu defines a formal trace function $\widehat{{\rm Tr}}_{W^\ell}(\ast:\tau)$ on $W^\ell$ by 
\begin{equation}
  \widehat{{\rm Tr}}_{W^\ell}(v:\tau):={\rm Tr}_{W^\ell} o(v) e^{2\pi i\tau (L(0)-c/24)},  \label{Intro1}
\end{equation}
and proves that these functions are well-defined as analytic functions on the upper half-plane 
$\CH=\{\tau\in \C\mid {\rm Im}(\tau)>0\}$. He then shows that for each $\gamma = \left( \begin{smallmatrix} a&b\\c&d \end{smallmatrix}\right) \in \SL_2(\Z)$ there exists a complex matrix $(A^\gamma_{ij})$ such that 
\begin{equation}
  \widehat{{\rm Tr}}_{W^\ell}\left(v:\frac{a\tau +b}{c\tau +d}\right)=(c\tau +d)^{\wt[v]}\sum_{k=1}^r A^\gamma_{\ell k}\widehat{{\rm Tr}}_{W^k}(v:\tau) \label{ZhuThm}
\end{equation}
for any $\tau\in \CH$ and $\wt[\cdot]$-homogeneous element $v\in V$.  
Here $\wt[\cdot]$ is the weight given by the Virasoro element $\widetilde{\omega}$ of the coordinate transformation VOA structure
$(V,Y[\cdot],\1,\widetilde{\omega})$ on $V$, which is given by setting
\begin{align*}
Y[v,z]&=Y(v,e^{2\pi iz}-1)e^{2\pi iz\wt(v)}=\sum_{n\in \mathbb{Z}} v[n]z^{-n-1}, \mbox{ and} \\
\widetilde{\omega}&=(2\pi i)^2 \left(\omega-\frac{c}{24}\1\right)\in V_{[2]}.
\end{align*}

We note that the action of $\SL_2(\Z)$ on $\CH$ is generated by an inversion 
$S=\left(\begin{smallmatrix} 0&-1\\1&0 \end{smallmatrix}\right):\tau \to \frac{-1}{\tau}$ and a parallel translation $T=\left(\begin{smallmatrix} 1&1\\0&1 \end{smallmatrix}\right):\tau\to \tau+1$. The invariance property of trace functions for $T$ follows easily from the structure of simple modules. Meanwhile, the matrix $(A^S_{ij})$ in (\ref{ZhuThm}) produced by the matrix $S$ contains interesting and exploitable information about $V$ and its modules. It is often called the $S$-matrix of $V$ and is denoted by $(s_{ij})$ rather than $(A^S_{ij})$.

As alluded to above, our motivation stems from observing that the power of $e$ in (1.1), i.e. 
$2\pi i\tau(L(0)-c/24)$, is the grade-preserving operator $o(\tau \widetilde{\omega}/2\pi i)$ 
of an element in a one-dimensional Jordan subalgebra $\C \widetilde{\omega}$ of $V_{[2]}$. 
We therefore treat cases where $V_{[2]}$ (and also $V_2$) contain a larger Jordan subalgebra $\CG$. Then, for $u\in \CG$ and $v\in V$, we define a multivariable formal trace function ${\rm Tr}_{W^\ell}o(v)e^{o(u)/2\pi i}$ and establish a new $\SL_2(\Z)$-invariance property for these functions. 

The first case we consider is that of an associative Jordan subalgebra. 
Let $V=(V,Y,\1,\omega)$ be a regular VOA and $\omega=e^1+\cdots+e^g$, 
where $e^j$ are conformal vectors that are mutually orthogonal with respect to the bilinear form discussed above. Here, an element $e\in V_2$ is called a conformal vector if $e$ is a Virasoro element of 
the subVOA generated by $e$, which we denote ${\rm VOA}(e)$. 
In this case, $\oplus_{j=1}^g\C e^j$ is an associative Jordan subalgebra of $V_2$. 
Set $\widetilde{e}^j =(2\pi i)^2 \left(e^j-\frac{c_j}{24}\1\right)$, where $c_j$ is the central charge of $e^j$.
Under this setting, for a grade-preserving operator $\alpha$ and 
$(\tau_1,\ldots,\tau_g)\in \CH^{g}$, we define a multivariable formal trace function by 
\begin{equation}
 \widehat{{\rm Tr}}_{W^h}(\alpha:\tau_1,\ldots,\tau_g):={\rm Tr}_{W^h}
\alpha e^{o\left(\sum_{j=1}^g \tau_j\widetilde{e}^j/2\pi i\right)}. \label{MultiTraceFct1}
\end{equation}
 If an element $u\in V$ is homogeneous with respect to the grading induced by the operator $L_j[0]:=\widetilde{e}^j[1]$, we denote its weight under this operator by $\wt_j[u]$. We say an element of $V$ is multi-$\prod\wt_j[]$-homogeneous if it is homogeneous with respect to $L_j[0]$ for all $j$. For a multi-$\prod\wt_j[]$-homogeneous element $w\in V$, let $\otimes_{j=1}^g {\rm VOA}(e^j)w$ denote the $\otimes_{j=1}^g{\rm VOA}(e^j)$-submodule generated by $w$. Then we have the following theorem, which is proved in Section \ref{Section-Preliminary}.

\begin{thm}\label{MMT}
 Let $V$ be a regular VOA and $\omega=\sum_{j=1}^ge^j$ be a 
decomposition of the Virasoro element $\omega$ by mutually orthogonal conformal vectors $e^j$. Let $w\in V$ be a multi-$\prod\wt_j[]$-homogeneous element and assume the functions $\widehat{{\rm Tr}}_{W^h}\left(o(v):\tau_1,\ldots,\tau_g\right)$ are well-defined as analytic functions on $\CH^{g}$ for $v\in \otimes_{j=1}^g {\rm VOA}(e^j)w$. Then
$$ \widehat{{\rm Tr}}_{W^\ell}\left(o(v):\frac{a\tau_1 +b}{c\tau_1 +d},\ldots,\frac{a\tau_g +b}{c\tau_g +d}\right)=
\prod_{p=1}^g(c\tau_p +d)^{\wt_p[v]}\sum_{h=1}^r A^\gamma_{\ell h}\widehat{{\rm Tr}}_{W^h}\left(o(v):\tau_1,\ldots,\tau_g\right) $$
for $(\tau_1,\ldots,\tau_g)\in \CH^{g}$, where $(A^\gamma_{ij})$ is the matrix given in (\ref{ZhuThm}). Additionally, if $e^j(n)w=0$ for $n>1$, $1\leq j\leq g$, and the functions $\widehat{{\rm Tr}}_{W^h}\left(o(w):\tau_1,\ldots,\tau_g\right)$ are well-defined as analytic functions on $\CH^{g}$, then so are the functions $\widehat{{\rm Tr}}_{W^h}\left(o(v):\tau_1,\ldots,\tau_g\right)$ for any $v\in \otimes_{j=1}^g {\rm VOA}(e^j)w$.
\end{thm}

One important instance when $\omega$ decomposes as stated in the previous theorem is when $V$ contains a simple, regular subVOA $U=(U,Y,\1, e)$, and we additionally consider the commutant of $U$ in $V$ given by $U^c=(U^c={\rm Com}_V(U),Y,\1 ,f =\omega -e)$, where ${\rm Com}_V(U):=\{v\in V \mid u_nv=0\mbox{ for }u\in U, n\in\N\}$. As an application of Theorem \ref{MMT}, we prove the following theorem in Section \ref{Section-Commutant}.

\begin{thm}\label{MC} 
Let $V$ be a simple regular VOA and $U$ a simple regular subVOA of $V$. Suppose also that the commutant $U^c$ of $U$ is simple, regular, and satisfies $(U^c)^c =U$, and that $U$, $U^c$, and $V$ are all of CFT-type and self-dual. Then all simple $U$-modules appear in some simple $V$-module.
\end{thm}

In the second case, we consider a Jordan algebra of type $B_g$. That is, 
a Jordan algebra isomorphic to the space ${\frak{S}}_g(\C)$ consisting of all symmetric complex matrices of degree $g$. More specifically, we have $V_2$ contains a Griess subalgebra $\CG:=\oplus_{1\leq i\leq j\leq g} \C \omega^{ij}$ and there exists an algebra isomorphism $\mu:{\frak{S}}_g(\C)\to \CG$ 
satisfying $\mu(E_{ij}+E_{ji})=2\omega^{ij}$ ($=\omega^{ij}+\omega^{ji}$) and $\mu(I_g)=\omega$,  
where $E_{ij}$ denotes an elementary matrix which has $1$ in the $(i,j)$-entry and zeros elsewhere. 
Here we call a subalgebra $\CG$ of $(V_2,\times_1)$ a Griess subalgebra if $v_2u=0$ for $v,u\in \CG$, 
where a $1$-product $u\times_1 v$ is given by $u_1v$. Such a subalgebra is commutative, but not necessarily associative. We note that the product $(V_2,\times_1)$ is a $\mathbb{C}$-algebra. If $V_0=\mathbb{C}\textbf{1}$ and $V_1=0$, then $V_2$ becomes a Griess (sub)algebra of $V$ (see, for example, \cite{AM}). Such Griess algebra structures are generalizations of the original Griess algebra, which is the commutative non-associative algebra on a real vector space of dimension $19886$ whose automorphism group is the Monster group. Trace functions of many variables associated with Griess algebras additionally appear in \cite{M2}. We also mention that in our setting we have $u_3 v$ is the bilinear form discussed above, that is, $u_3 v=\langle u,v\rangle \textbf{1}$.

In this situation, for $A=(\tau_{ij})\in \CH_g$ we define a multivariable trace function   
\begin{equation}
 \widehat{{\rm Tr}}_{W^\ell}\left(o(v):A\right)
={\rm Tr}_{W^\ell}o(v)e^{o\left(2\pi i \left(\mu(A)-\frac{{\rm tr}(A)c}{24g} \right) \right)}, \label{MultiTraceFct2}
\end{equation}
where $\CH_g=\{X+Yi\mid X,Y\in {\frak{S}}_g(\R), Y \mbox{ is positive definite}\}$ 
is the Siegel upper half-space. The action of $\SL_2(\Z)=\langle T,S \rangle$ on $\CH_g$ is given by 
$T(Z)=Z+E_g$ and $S(Z)=-Z^{-1}$ for $Z\in \CH_g$, where $E_g$ is the $g\times g$ identity matrix. 
Our next result, which is found in Section \ref{Section-Jordan}, establishes the invariance for a Siegel-type inversion.

\begin{thm}\label{SMT2}
Suppose $\widehat{{\rm Tr}}_{W^j}\left(o(\1):A\right)$ is a well-defined analytic function on $\CH_g$ for $j=1,\ldots,r$. Then 
\begin{equation}
\widehat{{\rm Tr}}_{W^j}\left(o(\1):-A^{-1}\right)=\sum_{h=1}^r s_{jh}\widehat{{\rm Tr}}_{W^h}\left(o(\1):A\right), \notag
\end{equation}
where $(s_{jh})$ is the $S$-matrix given for $\gamma =S$ in (\ref{ZhuThm}). \hfill \qed
\end{thm}

There are many known VOAs containing a Jordan algebra of type $B_g$. 
For example, a VOA $M(1)^{\otimes g}$ of free boson type constructed 
from a $g$-dimensional vector space $\C^g$ and its fixed point subVOA 
$(M(1)^{\otimes g})^+$ by an automorphism $-1$ on $\C^g$ 
contain a Griess subalgebra isomorphic to a Jordan algebra of type $B_g$ 
(that is, a Jordan algebra consisting of all symmetric complex matrices of degree $g$). 
The famous moonshine VOA $V^{\natural}$ also contains a Griess subalgebra isomorphic to 
a Jordan algebra of type $B_{24}$. Moreover, $V^{\natural}$ has only one simple module and its $S$-matrix is $(s_{ij})=I_1$.
We also note that the second author and Ashihara have shown in \cite{AM} that for any $c\in \C$ and $g\in \N$, 
there is a VOA $AM(g,c)$ with central charge $c$ whose Griess algebra is a Jordan algebra of 
type $B_g$. 

We conclude this paper with Section \ref{Section-Applications}, where we apply the above results to prove the inversion transformation property and convergence for ordinary Siegel theta series. See Proposition \ref{PropLattice} below for a detailed statement of this result.

\section{Preliminaries and simultaneous transformations} \label{Section-Preliminary}
We first recall the following notation and results from \cite{Z}. 
Since we will treat power series of $e^{2\pi i\tau_j}$ for various $\tau_j\in \CH$, 
we denote the $q$-power expansion of Eisenstein series $G_{2k}(\tau)$ 
by $\widetilde{G}_{2k}(\tau)$, where $q=e^{2\pi i\tau}$. 
Namely, 
\begin{equation}
\widetilde{G}_{2k}(\tau)=2\zeta(2k)+\frac{2(2\pi i)^{2k}}{(2k-1)!}\sum_{n=1}^{\infty}\sigma_{2k-1}(n)q^n. 
\label{EisensteinExpansion} %2.1
\end{equation}
Under the the action of a matrix $\gamma = \left(\begin{smallmatrix} a&b\\c&d \end{smallmatrix}\right) \in \text{SL}_2 (\mathbb{Z})$, these transform as 
\begin{equation}
 \begin{aligned}
\widetilde{G}_2 \left(\tfrac{a\tau +b}{c\tau +d} \right) &=(c\tau +d)^2\widetilde{G}_2(\tau)-2\pi i c(c\tau +d) \hspace{5mm} \text{  and  }\\
\widetilde{G}_{2k}\left(\tfrac{a\tau +b}{c\tau +d}\right)&=(c\tau +d)^{2k}\widetilde{G}_{2k}(\tau) \hspace{5mm} \text{  for  }k>1.
\end{aligned} \label{EisensteinTransforms} %2.2
\end{equation}

One of the most important results in \cite{Z} is the following, which we will often use.  

\begin{lmm}\label{Zhu}
For any VOA $V$ and any $L(0)$-gradable module $M$ whose grading is bounded below, we have
\begin{align}
\widehat{{\rm Tr}}_M\left(o(a[0]b):\tau\right)&= 0, \hspace{5mm} \text{  and  } \label{ZhuRec0} \\
\widehat{{\rm Tr}}_M\left(o(a)o(b):\tau\right)&=\widehat{{\rm Tr}}_M\left(o(a[-1]b):\tau\right)
-\sum_{k=1}^{\infty}\widetilde{G}_{2k}(\tau)\widehat{{\rm Tr}}_M\left(o(a[2k-1]b):\tau\right). \label{ZhuRec1} %2.3
\end{align}
as formal complex power series of $e^{2\pi i\tau}$ for $a,b\in V$. 
\end{lmm}

In this section, we let $V=(V,Y(\ast,z),\1,\omega)$ be a regular VOA of CFT-type 
and assume that $V$ has 
a set $\{e^1,\ldots,e^g\}$ of mutually orthogonal conformal vectors such that 
$\omega=e^1+\cdots+e^g$. Then $\widetilde{e}^j$ 
are mutually orthogonal conformal vectors of the coordinate transformation 
VOA $(V,Y[\ast,z],\1,\widetilde{\omega})$ and $\widetilde{\omega}=\sum_{j=1}^g \widetilde{e}^j$, 
where $c_j$ is the corresponding central charges of $e^j$ and $\widetilde{e}^j=(2\pi i)^2 \left(e^j-\frac{c_j}{24}\1\right)$. 
Let $M$ be a $V$-module and recall the multivariable functions (\ref{MultiTraceFct1}).  
Clearly, we have 
\begin{equation}
\frac{\pd}{\pd \tau_j}\widehat{{\rm Tr}}_M\left(o(v):\tau_1,\ldots,\tau_g\right)
=\frac{1}{2\pi i}\widehat{{\rm Tr}}_M\left(o(\widetilde{e}^j)o(v):\tau_1,\ldots,\tau_g\right). \label{Partials1}
\end{equation}
Since $e^j$ are mutually orthogonal, $[o(\widetilde{e}^j),o(\widetilde{e}^h)]=0$ and so 
we have the commutativity of partial differentials,
\begin{equation}
\frac{\pd}{\pd \tau_j}\frac{\pd}{\pd \tau_h}\widehat{{\rm Tr}}_M\left(o(v):\tau_1,\ldots,\tau_g\right)
=\frac{\pd}{\pd \tau_h}\frac{\pd}{\pd \tau_j}\widehat{{\rm Tr}}_M\left(o(v):\tau_1,\ldots,\tau_g\right), \label{Partials2}
\end{equation} 
for any $j$ and $h$. We also note that 
$$\lim_{\forall\tau_i\to \tau}\widehat{{\rm Tr}}_{M}\left(o(v):\tau_1,\ldots,\tau_g\right)
=\widehat{{\rm Tr}}_{M}(v:\tau).$$

Using arguments as in \cite{Miy-Int}, we obtain the following result.

\begin{lmm}
 We have
\begin{equation}
\widehat{{\rm Tr}}_M\left(o(e^j[0]b):\tau_1,\ldots,\tau_g\right)=0,
\label{ZhuRec0a}
\end{equation}
and
\begin{equation}
\begin{aligned}
\widehat{{\rm Tr}}_M\left(o(e^j)o(b):\tau_1,\ldots,\tau_g\right)&=
\widehat{{\rm Tr}}_M\left(o(e^j[-1]b):\tau_1,\ldots,\tau_g\right) \\
&\hspace{10mm} -\sum_{k=1}^{\infty}\widetilde{G}_{2k}(\tau_j)
\widehat{{\rm Tr}}_M\left(o(e^j[2k-1]b):\tau_1,\ldots,\tau_g\right).
\end{aligned}
\label{ZhuRec3}
\end{equation}
\end{lmm}

\pr  
 For any $b\in V$, $o_k(b):=b_{\wt (b) -1-k}$ and $b_m$ are given by the vertex operator $Y^M(b,z)$ of $V$ on $M$.  
We note that all of the forthcoming actions are given by $Y^M$ of $V$ on $M$. For formal variables $z_1$ and $z_2$ we consider functions $F_M$ defined by
\[
F_M \left((v^1 ,z_1),(v^2, z_2),q_1,\dots ,q_g \right)
= \Tr_M z_1^{\wt(v^1)}z_2^{\wt(v^2)}Y^M(v^1,z_1)Y^M(v^2,z_2)q_1^{o(e^1)}\dots q_g^{o(e^g)},
\]  
where $q_j=e^{2\pi i\tau_j}$. See \cite{Miy-Int} for more about such functions as well as many of the ideas we will use here (though note that our notation $o_k (b)$ equals $o_{-k}(b)$ there). In fact, since $o(e^1)$ commutes with each $o(e^r)$ for $1\leq r\leq g$, we obtain (\ref{ZhuRec0a}) just as in Proposition 3.1 of \cite{Miy-Int}.

The key point in establishing (\ref{ZhuRec3}) is that $o_k(e^1)$ commutes with $o(e^r)$ for $r\geq 2$ and $[o_k(e^1),o(e^1)]=-k o_k(e^1)$. Recall our Virasoro notation, $L_j (n):= e^j_{n+1}$, which implies $o_k(e^1)=L_1(-k)$. Therefore, we have
\begin{align*}
  &\Tr_M o_k(e^1)o_{-k}(b)q_1^{o(e^1)}\dots q_g^{o(e^g)} \\
&=\Tr_M [o_k(e^1),o_{-k}(b)] q_1^{o(e^1)}\dots q_g^{o(e^g)} 
 +\Tr_M o_{-k}(b)o_k(e^1)q_1^{o(e^1)}\dots q_g^{o(e^g)} \\
&=\Tr_M [e^1_{1-k},b_{\wt (b)-1+k}] q_1^{o(e^1)}\dots q_g^{o(e^g)} 
 +\Tr_M o_{-k}(b)q_1^{o(e^1)-k}q_2^{o(e^2)}\dots q_g^{o(e^g)}o_k(e^1) \\
&=\Tr_M \left(\sum_{j=0}^{\infty} \binom{1-k}{j}o\left(e^1_jb\right) \right) q_1^{o(e^1)}\dots q_g^{o(e^g)} 
 +\Tr_M o_k(e^1)o_{-k}(b)q_1^{o(e^1)-k}q_2^{o(e^2)}\dots q_g^{o(e^g)} .
 \end{align*}
 Rearranging, we find
\[
\left(1-q_1^{-k}\right)\Tr_M o_k(e^1)o_{-k}(b)q_1^{o(e^1)}\dots q_g^{o(e^g)} 
=\Tr_M \left(\sum_{j=0}^{\infty} \binom{1-k}{j}o\left(e^1_jb\right)\right) q_1^{o(e^1)}\dots q_g^{o(e^g)}, 
\]
and for $k\not =0$,
\begin{equation}
\begin{aligned}
\Tr_M &\left\{o_k(e^1)o_{-k}(b)q_1^{o(e^1)}\right\}q_2^{o(e^2)} \dots q_g^{o(e^g)} \\
&=\Tr_M \left\{ \frac{1}{1-q_1^{-k}}\left(\sum_{j=0}^{\infty} \binom{1-k}{j}o\left(e^1_jb\right) \right)q_1^{o(e^1)} \right\} q_2^{o(e^2)}\dots q_g^{o(e^g)} .
\end{aligned}\label{new1}
\end{equation}
Again, the important fact is that $o_k(e^1)$ commutes with $o(e^r)$ ($r\not=1$) so that only $q_1^{-k}$ arises. 
We note that the expression contained in $\{\dots \}$ of (\ref{new1}) is the same as in Proposition 3.2 of \cite{Miy-Int} by viewing $q_1^{o(e^1)}$ as $q^{L(0)}$. Therefore, by the same argument found there, we have 
\begin{align*}
&F_M((e^1,x),(b,z),q_1,\dots ,q_g) \\
&=\Tr_M o(e^1)o(b)q_1^{o(e^1)}\dots q_g^{o(e^g)} 
+\sum_{m\in \N}P_{m+1}\left(\frac{z}{x},q \right) \Tr_M o(e^1[m]b)q_1^{o(e^1)}\dots q_g^{o(e^g)} 
\end{align*}
and
\begin{align*}
&F_M ((b,z),(e^1,x),q_1,\dots ,q_g) \\
&=\Tr_M o(e^1)o(b)q_1^{o(e^1)}\dots q_g^{o(e^g)} 
+\sum_{m\in \N}\left\{P_{m+1}\left(\frac{zq}{x},q \right)-\delta_{m,0} \right\}\Tr_M o(e^1[m]b)q_1^{o(e^1)}\dots q_g^{o(e^g)} ,
\end{align*}
where the functions $P_k (z,\tau)$ are defined in \cite{Miy-Int} and are also the similarly denoted functions in \cite{Z} multiplied by $(2\pi i)^{-k}$. By using the associativity of endomorphisms for $Y^M$, we find
\begin{align*}
    o(e^1_kb)&=(e^1_kb)_{\wt(b)-k}\\
&=\sum_{j=0}^{\infty}\binom{k}{j}(-1)^j \left\{ e^1_{k-j}b_{\wt(b)-k+j}-(-1)^kb_{\wt(b)-j}e^1_j \right\} \\
&=\sum_{j=0}^{\infty}\binom{k}{j}(-1)^j \left\{ o_{-k+j+1}(e^1)o_{k-j-1}(b)-(-1)^ko_{j-1}(b)o_{-j+1}(e^1)\right\}.  
\end{align*}
Replacing $q^{L(0)}$ in Equation (3.7) of Proposition 3.3 in \cite{Miy-Int} by $q_1^{o(e^1)}$ and  multiplying by $q_2^{o(e^2)}\dots q_g^{o(e^g)}$, we obtain a similar result. Namely, 
\begin{align*}
&\Tr_M o(e^1)o(b)q_1^{o(e^1)}\dots q_g^{o(e^g)} \\
&=\Tr_M o(e^1[-1]b)q_1^{o(e^1)}\dots q_g^{o(e^g)}
-\sum_{k=1}^{\infty}\widetilde{G}_{2k}(\tau_1)\Tr_M o(e^1[2k-1]b)q_1^{o(e^1)} \dots q_g^{o(e^g)},
\end{align*}
giving (\ref{ZhuRec3}). \hfill \qed  

We are now in position to prove Theorem \ref{MMT}.

%*********
%*********Proof of Theorem 1
%*********
\par \vspace{3mm}\noindent [{\bf Proof of Theorem \ref{MMT}}] \hspace{3mm}
In this proof, $L_k[m]$ denotes $\widetilde{e}^k[m+1]$.  
To simplify the notation, we will write the proof for the case $g=2$, 
but there is no difference for $g\geq 3$. \\
\indent We first prove the statement that $\widehat{{\rm Tr}}_{W^h}\left(o(v):\tau_1,\tau_2\right)$ is a well-defined analytic function on $\CH^{2}$, so long as $\widehat{{\rm Tr}}_{W^h}\left(o(w):\tau_1,\ldots,\tau_g\right)$ is, and $L_j(n)w=0$ for $n\geq 1$, $j=1,2$. We do so by induction on $\wt_1[v]+\wt_2[v]$, after assuming this is true for the base case $v=w$. More generally, any $v\in \otimes_{j=1}^2 {\rm VOA}(e^j)w$ is of the form $v= \otimes_{j=1}^2 L_j [-m_{1_j}] \cdots L_j [-m_{d_j}] w$ for $m_{i_j} \geq 1$. Since $L_j [-n]$ is generated by $L_j [-1]$ and $L_j [-2]$, we 
may take $m_{i_j}=1,2$. 
Moreover, by (\ref{ZhuRec0a}) we may assume $m_{1_j}=2$. Since $\widehat{{\rm Tr}}_{W^h}\left(o(\otimes_{j=1}^2 \widehat{L_j [-m_{1_j}]} L_j [-m_{2_j}] \cdots L_j [-m_{d_j}] w):\tau_1,\tau_2\right)$, where $\widehat{L_j [-m_{1_j}]}$ denotes the omission of one or both of these terms, is analytic by our induction hypothesis, then (\ref{Partials1}) and (\ref{ZhuRec3}) imply $\widehat{{\rm Tr}}_{W^h}\left(o(v):\tau_1,\tau_2\right)$ is also analytic.\\
\indent We now turn to proving the functional equation. Set $\gamma =\left( \begin{smallmatrix} a&b \\ c&d \end{smallmatrix}\right) \in \text{SL}_2 (\mathbb{Z})$. We will prove 
\begin{equation}
\begin{aligned}
 \widehat{{\rm Tr}}_{W^\ell}&\left(o(v):\frac{a\tau_1 +b}{c\tau_1 +d},\frac{a\tau_2 +b}{c\tau_2 +d} \right) \\
 &=(c\tau_1 +d)^{\wt_1[v]}(c\tau_2 +d)^{\wt_2[v]} \sum_{h=1}^r A^\gamma_{\ell h} \widehat{{\rm Tr}}_{W^h}\left(o(v):\tau_1,\tau_2\right).
\end{aligned} \label{MMTproof1}
\end{equation}
For ease of notation, we set $\gamma \tau_i := \frac{a\tau_i +b}{c\tau_i +d}$ and $j(\gamma, \tau_i):=(c\tau_i +d)$. To begin, we consider $(\tau_1,\tau_2)=(\tau,\tau)$ as a base point. 
Since $\lim_{\forall \tau_i\to \tau}\widehat{{\rm Tr}}_{W^\ell}\left(o(v):\tau_1,\tau_2\right)=\widehat{{\rm Tr}}_{W^\ell} \left(v:\tau\right)$, 
Zhu's theorem (cf.\ (\ref{ZhuThm})) implies  
\begin{equation}
\begin{aligned}
  \lim_{\forall \tau_i\to \tau}\widehat{{\rm Tr}}_{W^\ell} & \left(o(v):\gamma \tau_1,\gamma \tau_2\right)\\
&=\sum_{h=1}^r A^\gamma_{\ell h}\lim_{\forall \tau_i\to \tau}j(\gamma , \tau_1)^{\wt_1[v]}j(\gamma , \tau_2)^{\wt_2[v]}\widehat{{\rm Tr}}_{W^h}\left(o(v):\tau_1,\tau_2\right).
\end{aligned} \label{MMTproof2}
\end{equation}
Namely, (\ref{MMTproof1}) is true for the base point $(\tau,\tau)$. 
We will next show that higher order partial derivatives on both sides of (\ref{MMTproof2}) by 
$\tau_1$ and $\tau_2$ still coincide with each other when evaluated at $(\tau,\tau)$. This in turn implies the Taylor series expansions about $(\tau ,\tau)$ of the analytic left and right hand sides of (\ref{MMTproof1}) are equal on a neighborhood about $(\tau ,\tau)$, and thus on all of $\CH^{2}$.\\ 
\indent To simplify the arguments, we will prove the equality for higher order partial 
derivatives by $\tau_1$. Namely, we will prove 
\begin{equation}
\begin{aligned}
  \lim_{\forall \tau_i\to \tau} &\frac{\pd^p}{\pd\tau_1^p} \left[
\widehat{{\rm Tr}}_{W^\ell}\left(o(v):\gamma \tau_1,\gamma \tau_2\right)\right] \\
& =\lim_{\forall \tau_i\to \tau}\frac{\pd^p}{\pd\tau_1^p}\left[ j(\gamma , \tau_1)^{\wt_1[v]}j(\gamma , \tau_2)^{\wt_2[v]} \sum_{h=1}^r A^\gamma_{\ell h}\widehat{{\rm Tr}}_{W^h}\left(o(v),\tau_1,\tau_2\right)\right]
\end{aligned}
\label{MMTproof3}
\end{equation}
for any $p\in \N$ and $v\in \otimes_{j=1}^2{\rm VOA}(e^j)w$ by induction. 
For the combinations with $\frac{\pd}{\pd \tau_2}$, we can prove the assertion by using (\ref{Partials1}),(\ref{Partials2}), and (\ref{ZhuRec3}). 
We note that (\ref{MMTproof3}) is true for $p=0$, and we next assume that it holds for all $v\in \otimes_{j=1}^2{\rm VOA}(e^j)w$ 
and $p \leq m$. In particular, we have 
\begin{equation}
\begin{aligned}
 \lim_{\forall \tau_i\to \tau} &\frac{\pd^m}{\pd\tau_1^m}\left[ j(\gamma , \tau_1)^{\wt_1[v]+2}j(\gamma , \tau_2)^{\wt_2[v]}
\sum_{h=1}^r A^\gamma_{\ell h}\widehat{{\rm Tr}}_{W^h}\left(o(L_1[-2]v):\tau_1,\tau_2\right)\right] \\
 &= \lim_{\forall \tau_i\to \tau}\frac{\pd^m}{\pd\tau_1^m}\left[
\widehat{{\rm Tr}}_{W^\ell}\left(o(L_1[-2]v):\gamma \tau_1,\gamma \tau_2\right)\right],
\end{aligned}
\label{MMTproof4}
\end{equation}
and more generally 
\begin{equation}
\begin{aligned}
 \lim_{\forall \tau_i\to \tau} &\frac{\pd^m}{\pd\tau_1^m}\left[ j(\gamma , \tau_1)^{\wt_1[v]+2-2k}j(\gamma , \tau_2)^{\wt_2[v]}
\sum_{h=1}^r A^\gamma_{\ell h}\widehat{{\rm Tr}}_{W^h}\left(o(L_1[2k-2]v):\tau_1,\tau_2\right)\right] \\
 &= \lim_{\forall \tau_i\to \tau}\frac{\pd^m}{\pd\tau_1^m}\left[
\widehat{{\rm Tr}}_{W^\ell}\left(o(L_1[2k-2]v):\gamma \tau_1,\gamma \tau_2\right)\right]
\end{aligned}
\label{MMTproof4b}
\end{equation}
for any $k\geq 0$. Furthermore, applying the Leibniz rule for higher order derivations on the product of $(c\tau_1 +d)^n \widetilde{G}_{2k}(\tau_1)$ with the left and right hand sides of (\ref{MMTproof4b}), we also have
\begin{equation}
\begin{aligned}
 \lim_{\forall \tau_i\to \tau} & \frac{\pd^m}{\pd\tau_1^m}\biggl [(c\tau_1 +d)^n \widetilde{G}_{2k}(\tau_1) j(\gamma , \tau_1)^{\wt_1[v]+2-2k}j(\gamma , \tau_2)^{\wt_2[v]} \\
&\hspace{15mm}\times \sum_{h=1}^r A^\gamma_{\ell h}\widehat{{\rm Tr}}_{W^h}\left(o(L_1[2k-2]v):\tau_1,\tau_2\right)\biggl ] \\
 &= \lim_{\forall \tau_i\to \tau} \frac{\pd^m}{\pd\tau_1^m}\left[(c\tau_1 +d)^n \widetilde{G}_{2k}(\tau_1)
\widehat{{\rm Tr}}_{W^\ell}\left(o(L_1[2k-2]v):\gamma \tau_1,\gamma \tau_2\right)\right]
\end{aligned}
\label{MMTproof5}
\end{equation}
for any $k,n\geq 0$, where we set $\widetilde{G}_0 (\tau_1):=1$.

Since $L_1[-2]=\widetilde{e}^1[-1]$, (\ref{ZhuRec3}) and direct calculation
gives
\begin{equation}
\begin{aligned}
\mbox{\rm RHS of (\ref{MMTproof4})}=&{\displaystyle \lim_{\forall \tau_i\to \tau}\frac{\pd^m}{\pd\tau_1^m} \left[
\widehat{{\rm Tr}}_{W^\ell} \left(o(\widetilde{e}^1)o(v):\gamma \tau_1,\gamma \tau_2\right) \right]}\cr
&\hspace{0mm}+\lim_{\forall \tau_i\to \tau}\frac{\pd^m}{\pd\tau_1^m}  \left[
\widetilde{G}_2(\gamma \tau_1)\wt_1[v]\widehat{{\rm Tr}}_{W^\ell} \left(o(v):\gamma \tau_1,\gamma \tau_2\right)\right]\\
&\hspace{0mm}+\lim_{\forall \tau_i\to \tau}\frac{\pd^m}{\pd\tau_1^m} \left[\sum_{k=2}^{\infty} \widetilde{G}_{2k}(\gamma \tau_1)
\widehat{{\rm Tr}}_{W^\ell} \left(o(L_1 [2k-2]v):\gamma \tau_1,\gamma \tau_2\right)\right].
\end{aligned}
\notag %\label{MMTproof6}
\end{equation}
Furthermore, from (\ref{Partials1}) we have
\begin{equation}
\begin{aligned}
 \frac{\pd^m}{\pd\tau_1^m}& \left[\widehat{{\rm Tr}}_{W^\ell} \left(o(\widetilde{e}^1)o(v):\gamma \tau_1,\gamma \tau_2 \right) \right]
= \frac{\pd^m}{\pd\tau_1^m}  2\pi i\left[ \frac{\pd}{\pd(\gamma \tau_1)}\widehat{{\rm Tr}}_{W^\ell} \left(o(v):\gamma \tau_1,\gamma \tau_2\right)\right] \\
&= 2\pi i\frac{\pd^m}{\pd\tau_1^m} \left[ j(\gamma , \tau_1)^2 \left[\frac{\pd}{\pd\tau_1}\widehat{{\rm Tr}}_{W^\ell} \left(o(v):\gamma \tau_1,\gamma \tau_2\right)\right]\right]\\
&=2\pi i\frac{\pd^{m+1}}{\pd\tau_1^{m+1}}\left[ j(\gamma , \tau_1)^2 \widehat{{\rm Tr}}_{W^\ell} \left(o(v):\gamma \tau_1,\gamma \tau_2\right)\right]\\
&\hspace{10mm}-\frac{\pd^m}{\pd\tau_1^m}\left[ 4\pi icj(\gamma , \tau_1) \widehat{{\rm Tr}}_{W^\ell} \left(o(v):\gamma \tau_1,\gamma \tau_2\right)\right],
\end{aligned}
\notag %\label{MMTproof7}
\end{equation}
so that
\begin{equation}
\begin{aligned}
\mbox{\rm RHS of (\ref{MMTproof4})}&=\displaystyle \lim_{\forall \tau_i\to \tau}2\pi i\frac{\pd^{m+1}}{\pd\tau_1^{m+1}}\left[ j(\gamma , \tau_1)^2 \widehat{{\rm Tr}}_{W^\ell} \left(o(v):\gamma \tau_1,\gamma \tau_2\right)\right]\\
&\hspace{5mm}- \lim_{\forall \tau_i\to \tau}\frac{\pd^m}{\pd\tau_1^m}\left[ 4\pi icj(\gamma , \tau_1) \widehat{{\rm Tr}}_{W^\ell} \left(o(v):\gamma \tau_1,\gamma \tau_2\right)\right]\\
&\hspace{5mm}+\lim_{\forall \tau_i\to \tau}\frac{\pd^m}{\pd\tau_1^m}\left[  
\widetilde{G}_2(\gamma \tau_1)\wt_1[v]\widehat{{\rm Tr}}_{W^\ell} \left(o(v):\gamma \tau_1,\gamma \tau_2\right)\right]\\
&\hspace{5mm}+\lim_{\forall \tau_i\to \tau}\frac{\pd^m}{\pd\tau_1^m} \left[\sum_{k=2}^{\infty} \widetilde{G}_{2k}(\gamma \tau_1)
\widehat{{\rm Tr}}_{W^\ell} \left(o(L_1 [2k-2]v):\gamma \tau_1,\gamma \tau_2\right)\right].
\end{aligned}
\notag %\label{MMTproof6}
\end{equation}
On the other hand, by (\ref{Partials1}) and (\ref{ZhuRec3}) we find
$$\begin{array}{l}
\mbox{\rm LHS of (\ref{MMTproof4})} \cr
={\displaystyle \lim_{\forall \tau_i\to\tau} \frac{\pd^m}{\pd\tau_1^m} \left[
j(\gamma , \tau_1)^{\wt_1[v]+2}j(\gamma , \tau_2)^{\wt_2[v]}\sum_{h=1}^r A^\gamma_{\ell h}\widehat{{\rm Tr}}_{W^h}
\left(o(\widetilde{e}^1)o(v):\tau_1,\tau_2\right)\right]}\cr
{\displaystyle \quad +
\lim_{\forall \tau_i\to\tau} \frac{\pd^m}{\pd\tau_1^m}\left[j(\gamma , \tau_1)^{\wt_1[v]+2}j(\gamma , \tau_2)^{\wt_2[v]}
\widetilde{G}_2(\tau_1)\wt_1[v]
\sum_{h=1}^r A^\gamma_{\ell h}\widehat{{\rm Tr}}_{W^h}\left(o(v):\tau_1,\tau_2\right)\right] } \cr
\mbox{}\quad{\displaystyle +\lim_{\forall \tau_i\to \tau} \frac{\pd^m}{\pd\tau_1^m}\left[j(\gamma , \tau_1)^{\wt_1[v]+2}j(\gamma , \tau_2)^{\wt_2[v]}
\sum_{k=2}^{\infty} \widetilde{G}_{2k}(\tau_1)\sum_{h=1}^r A^\gamma_{\ell h}\widehat{{\rm Tr}}_{W^h}\left(o(L_1 [2k-2]v):\tau_1,\tau_2\right) \right]}
\end{array}$$
$$\begin{array}{l}
={\displaystyle \lim_{\forall \tau_i\to \tau}\frac{\pd^m}{\pd\tau_1^m}\left[ j(\gamma , \tau_1)^{\wt_1[v]+2}j(\gamma , \tau_2)^{\wt_2[v]}
2\pi i\frac{\pd}{\pd\tau_1}\left[\sum_{h=1}^r A^\gamma_{\ell h}\widehat{{\rm Tr}}_{W^h}\left(o(v):\tau_1,\tau_2\right)\right]\right]}\cr
\mbox{}\quad{\displaystyle+\lim_{\forall \tau_i\to\tau}\frac{\pd^m}{\pd\tau_1^m}\left[j(\gamma , \tau_1)^{2}\widetilde{G}_2(\tau_1)
j(\gamma , \tau_1)^{\wt_1[v]}j(\gamma , \tau_2)^{\wt_2[v]}\wt_1[v]\sum_{h=1}^r A^\gamma_{\ell h}\widehat{{\rm Tr}}_{W^h}\left(o(v):\tau_1,\tau_2\right)\right]} \cr
\mbox{}\quad{\displaystyle +\lim_{\forall \tau_i\to \tau} \frac{\pd^m}{\pd\tau_1^m}\biggl [
\sum_{k=2}^{\infty} j(\gamma , \tau_1)^{2k}\widetilde{G}_{2k}(\tau_1)j(\gamma , \tau_1)^{\wt_1[v]+2-2k}j(\gamma , \tau_2)^{\wt_2[v]} }\cr
\mbox{}\hspace{30mm}\quad{\displaystyle \times \sum_{h=1}^r A^\gamma_{\ell h}\widehat{{\rm Tr}}_{W^h}
\left(o(L_1[2k-2]v):\tau_1,\tau_2\right)\biggl ]} 
\end{array}$$
$$\begin{array}{l}
={\displaystyle 2\pi i \lim_{\forall \tau_i\to \tau}\frac{\pd^{m+1}}{\pd\tau_1^{m+1}}\left[
 j(\gamma , \tau_1)^{\wt_1[v]+2}j(\gamma , \tau_2)^{\wt_2[v]} \sum_{h=1}^r A^\gamma_{\ell h}\widehat{{\rm Tr}}_{W^h}\left(o(v):\tau_1,\tau_2\right)\right]}\cr
\mbox{}\quad{\displaystyle -2\pi i \lim_{\forall \tau_i\to \tau}\frac{\pd^{m}}{\pd\tau_1^{m}}\biggl [
c(\wt_1[v]+2)j(\gamma , \tau_1)^{\wt_1[v]+1}j(\gamma , \tau_2)^{\wt_2[v]}} \cr
\mbox{}\hspace{30mm}\quad{\displaystyle \times \sum_{h=1}^r A^\gamma_{\ell h}\widehat{{\rm Tr}}_{W^h}\left(o(v):\tau_1,\tau_2\right) \biggl ]}\cr
\mbox{}\quad{\displaystyle +\lim_{\forall \tau_i\to \tau} \frac{\pd^m}{\pd\tau_1^m}\biggl [\left( \widetilde{G}_2(\gamma \tau_1)+2\pi icj(\gamma , \tau_1)\right) \wt_1[v]j(\gamma , \tau_1)^{\wt_1[v]}j(\gamma , \tau_2)^{\wt_2[v]} } \cr
 \mbox{}\hspace{30mm} \quad{\displaystyle \times \sum_{h=1}^r  A^\gamma_{\ell h}\widehat{{\rm Tr}}_{W^h}\left(o(v):\tau_1,\tau_2\right)\biggl ]}\cr
\mbox{}\quad{\displaystyle +\lim_{\forall \tau_i\to \tau} \frac{\pd^m}{\pd\tau_1^m}\biggl [
\sum_{k=2}^{\infty}\widetilde{G}_{2k}(\gamma \tau_1) j(\gamma , \tau_1)^{\wt_1[v]+2-2k}j(\gamma , \tau_2)^{\wt_2[v]}
} \cr
\mbox{}\hspace{30mm}\quad{\displaystyle \times \sum_{h=1}^r A^\gamma_{\ell h}\widehat{{\rm Tr}}_{W^h}\left(o(L_1 [2k-2]v):\tau_1,\tau_2\right)\biggl ]} \cr
={\displaystyle 2\pi i \lim_{\forall \tau_i\to \tau}\frac{\pd^{m+1}}{\pd\tau_1^{m+1}}
\left[j(\gamma ,\tau_1)^2 j(\gamma , \tau_1)^{\wt_1[v]}j(\gamma , \tau_2)^{\wt_2[v]} \sum_{h=1}^r A^\gamma_{\ell h}\widehat{{\rm Tr}}_{W^h}\left(o(v):\tau_1,\tau_2\right)\right]}\cr
\mbox{}\quad{\displaystyle -\lim_{\forall \tau_i\to \tau}\frac{\pd^{m}}{\pd\tau_1^{m}}\left[
4\pi icj(\gamma , \tau_1)j(\gamma , \tau_1)^{\wt_1[v]}j(\gamma , \tau_2)^{\wt_2[v]}\sum_{h=1}^r A^\gamma_{\ell h}\widehat{{\rm Tr}}_{W^h}\left(o(v):\tau_1,\tau_2\right)\right]}\cr
\mbox{}\quad{{\displaystyle+\lim_{\forall \tau_i\to \tau} \frac{\pd^m}{\pd\tau_1^m}\left[ \widetilde{G}_2(\gamma \tau_1)\wt_1[v]j(\gamma , \tau_1)^{\wt_1[v]}j(\gamma , \tau_2)^{\wt_2[v]}
\sum_{h=1}^r A^\gamma_{\ell h}\widehat{{\rm Tr}}_{W^h}\left(o(v):\tau_1,\tau_2\right)\right]}} \cr
\mbox{}\quad{\displaystyle +\lim_{\forall \tau_i\to \tau} \frac{\pd^m}{\pd\tau_1^m}\biggl [
\sum_{k=2}^{\infty} \widetilde{G}_{2k}(\gamma \tau_1)
j(\gamma , \tau_1)^{\wt_1[v]+2-2k}j(\gamma , \tau_2)^{\wt_2[v]}} \cr
\mbox{}\hspace{30mm}\quad{{\displaystyle \times \sum_{h=1}^r A^\gamma_{\ell h}\widehat{{\rm Tr}}_{W^h}\left(o(L_1 [2k-2]v):\tau_1 ,\tau_2 \right)\biggl ],}} \end{array}$$ 
where we also used the transformations (\ref{EisensteinTransforms}).

Equating our calculations for the left and right hand sides of (\ref{MMTproof4}) and using (\ref{MMTproof5}), we obtain 
\begin{equation}
\begin{array}{l}
\displaystyle{\lim_{\forall \tau_i\to \tau}\frac{\pd^{m+1}}{\pd\tau_1^{m+1}}
\left[j(\gamma , \tau_1)^2 j(\gamma , \tau_1)^{\wt_1[v]}j(\gamma , \tau_2)^{\wt_2[v]}\sum_{h=1}^r A^\gamma_{\ell h}\widehat{{\rm Tr}}_{W^h}\left(o(v):\tau_1,\tau_2\right)\right]}\cr
\displaystyle{\mbox{}\qquad \qquad =
 \lim_{\forall \tau_i\to \tau}\frac{\pd^{m+1}}{\pd\tau_1^{m+1}} \left[j(\gamma , \tau_1)^2\widehat{{\rm Tr}}_{W^\ell} \left(o(v):\gamma \tau_1,\gamma \tau_2\right)\right]}. 
 \end{array}\label{mmtpf1}
 \end{equation}
Meanwhile, using our induction hypothesis (\ref{MMTproof3}) for $0\leq p\leq m$, together with the higher order Leibniz rule, we find
\begin{equation}
\begin{array}{l}
\displaystyle{\lim_{\forall \tau_i\to \tau}\frac{\pd^{m+1}}{\pd\tau_1^{m+1}}
\left[j(\gamma , \tau_1)^2 j(\gamma , \tau_1)^{\wt_1[v]}j(\gamma , \tau_2)^{\wt_2[v]}\sum_{h=1}^r A^\gamma_{\ell h}\widehat{{\rm Tr}}_{W^h}\left(o(v):\tau_1,\tau_2\right)\right]}\cr
\displaystyle{\mbox{}  =
 \lim_{\forall \tau_i\to \tau}\frac{\pd^{m+1}}{\pd\tau_1^{m+1}} \left[j(\gamma , \tau_1)^2\widehat{{\rm Tr}}_{W^\ell} \left(o(v):\gamma \tau_1,\gamma \tau_2\right)\right]}\cr
 \displaystyle{\mbox{} \qquad -
 \lim_{\forall \tau_i\to \tau}j(\gamma , \tau_1)^2 \frac{\pd^{m+1}}{\pd\tau_1^{m+1}} \left[\widehat{{\rm Tr}}_{W^\ell} \left(o(v):\gamma \tau_1,\gamma \tau_2\right)\right]}\cr
 \displaystyle{\mbox{} \qquad +\lim_{\forall \tau_i\to \tau}j(\gamma , \tau_1)^2 \frac{\pd^{m+1}}{\pd\tau_1^{m+1}}
\left[ j(\gamma , \tau_1)^{\wt_1[v]}j(\gamma , \tau_2)^{\wt_2[v]}\sum_{h=1}^r A^\gamma_{\ell h}\widehat{{\rm Tr}}_{W^h}\left(o(v):\tau_1,\tau_2\right)\right]}.
 \end{array} \notag
 \end{equation}
Applying (\ref{mmtpf1}) produces the desired result. This completes the proof of Theorem \ref{MMT}.   \hfill \qed

\section{Commutant decomposition} \label{Section-Commutant}

In this section we prove Theorem \ref{MC}. 

\par \vspace{3mm}\noindent [{\bf Proof of Theorem \ref{MC}}] \hspace{3mm}
Set $W=U^c$. Let $\{V=V^1 , \dots ,V^p \}$ denote the set of simple $V$-modules and $\{U=U^1,\dots ,U^g, \dots ,U^s \}$ and $\{W=W^1 ,\dots ,W^h ,\dots ,W^t \}$ be the sets of simple $U$-modules and simple $W$-modules, respectively, where the sets $\{U^1, \dots ,U^g \}$ and $\{W^1 ,\dots ,W^h \}$ denote the modules which appear in some simple $V$-module. We will prove the theorem by contradiction, and assume that $s>g$.\\
\indent Set $\sigma = \left( \begin{smallmatrix} 0&-1 \\1&0 \end{smallmatrix} \right)$, and let $S=(s_{ij})$, $S^U$, and $S^W$ denote the matrices $(A_{ij}^\sigma)$ given in (\ref{ZhuThm}) of the trace functions on $V$-modules, $U$-modules, and $W$-modules, respectively. Moreover, let $M_{r,s}(F)$ denote the set of $r\times s$-matrices with entries in $F$.\\
\indent Viewing $V^k$ as a $U\otimes W$-module, we have the existence of a matrix $R^k =(r_{ij}^k) \in M_{s,t}(\mathbb{N})$ such that $V^k$ decomposes as
\[
 V^k =\bigoplus_{i,j} r_{ij}^k (U^i \otimes W^j),
\]
where $r_{ij}^k (U^i \otimes W^j)$ denotes a direct product of $r_{ij}^k$ copies of $U^i \otimes W^j$. By the transformation property of Theorem \ref{MMT}, we have
\begin{equation}
 (S^U) R^k ({}^tS^W) =\sum_{i=1}^p s_{ki}R^i . \label{eqN1}
\end{equation}
 \indent For variables $x_1 ,\dots ,x_p$, let $\underline{x}$ denote the set of variables $\{x_1 ,\dots ,x_p\}$. Set $R(\underline{x})=R(x_1 ,\dots ,x_p) := \sum_{i=1}^p x_i R^i$, and view $R(\underline{x})$ as a matrix over $K:= \mathbb{C} (x_1 ,\dots ,x_p)$. Then replacing $R^k$ in (\ref{eqN1}) with $R(\underline{x})$, we find
 \begin{align}
 (S^U) R (\underline{x}) ({}^tS^W)&= \sum_{\ell =1}^p x_{\ell} \left [ (S^U) R^\ell ({}^tS^W) \right ]=\sum_{\ell =1}^p x_\ell \left( \sum_{i=1}^p s_{\ell i} R^i \right) =\sum_{i=1}^p \left(\sum_{\ell =1}^p s_{\ell i} x_\ell \right) R^i \notag \\
  &= \sum_{i=1}^p \widehat{x}_i R^i = R(\underline{\widehat{x}}), \label{eqN2}
%  &=\sum_{\ell =1}^p R (s_{\ell 1}x_1 ,\dots ,s_{\ell p} x_p) \notag \\
%  &=R \left( \sum_{\ell =1}^p s_{\ell 1}x_1 ,\dots ,\sum_{\ell =1}^p s_{\ell p}x_p  \right)=R(\widehat{x_1},\dots ,\widehat{x_p})=R\left(\underline{\widehat{x}}\right)
 \end{align}
 where we set $\widehat{x}_j =\sum_{\ell =1}^p s_{\ell j} x_\ell$. Moreover, we extend the transformation $(s_{ij}) \colon x_j \to \sum_{\ell =1}^p s_{\ell j}x_\ell =\widehat{x}_j$ to a $\mathbb{C}$-automorphism $\phi$ of $K$. We also let $\underline{x}^\phi$ denote the application of $\phi$ on each $x_i$. For example, 
 \[
 R(\underline{x}^\phi) =R\left( \phi (x_1) ,\dots ,\phi (x_p)\right) =R\left( \sum_{\ell =1}^p s_{\ell 1}x_\ell ,\dots ,\sum_{\ell =1}^p s_{\ell p}x_\ell \right).
 \]
  We can then rewrite (\ref{eqN2}) as
 \begin{equation}
 (S^U) R (\underline{x}) = R(\underline{x}^\phi)({}^tS^W)^{-1}. \label{eqN3}
\end{equation}
\indent It follows from the assumption $s>g$ that $r_{sj}^k =0$ for all $j,k$. Hence, the $s$-th row of $R(\underline{x})$ is zero, as is the $s$-th row of $R(\underline{x}^\phi)$. Therefore the $s$-th row of $R(\underline{x}^\phi)({}^tS^W)^{-1}$ on the right hand side of (\ref{eqN3}) is zero, and thus so is the $s$-th row of $(S^U) R (\underline{x})$. In particular, the $(s,1)$-entry of $(S^U) R (\underline{x})$ is zero. Explicitly, we have
\begin{equation}
 \sum_{j=1}^s S^U_{sj} R(\underline{x})_{j1} =0, \label{eqN11}
\end{equation}
where $S^U_{ij}$ and $R(\underline{x})_{ij}$ denote the $(i,j)$-entries of the matrices $S^U$ and $R(\underline{x})$, respectively. Meanwhile, $S^U_{s1}$ is nonzero by the Verlinde formula. Additionally, since $W=U^c$ and $U=W^c$, we have $r_{j1}^1 = 0$ for $j>1$. This implies $R(\underline{x})_{j1} \in \mathbb{N}[x_2,\dots ,x_p]$ for $j>1$. Finally, noting that $r^1_{11}=1$, we have $R(\underline{x})_{11} \in x_1 +\mathbb{N} [x_2 ,\dots ,x_p]$. It follows that 
\[
\sum_{j=1}^s S^U_{sj} R(\underline{x})_{j1} \in S^U_{s1}x_1 +\mathbb{C}[x_2 ,\dots ,x_p],
\]
which cannot equal zero. This contradicts (\ref{eqN11}), and the proof is complete. \hfill \qed

\section{Jordan algebra of type $B_g$} \label{Section-Jordan}
Let ${\frak{S}}_g(\C)$ denote the set of symmetric matrices of degree $g$ and 
$\CH_g$ be the Siegel upper half-space $\{X+Yi\mid X,Y\in {\frak{S}}_g(\R), Y \mbox{ is positive definite} \}$.   
We note that ${\frak{S}}_g(\C)$ is a Jordan algebra of type $B_g$. 

In this section, we assume that there is a Griess subalgebra 
$\CG\subseteq V_2$ isomorphic to a Jordan algebra ${\frak{S}}_g(\C)$ such that 
the identity matrix corresponds to $\omega$, and a primitive idempotent corresponds 
to a conformal element with a central charge $c/g$.  
We denote its ring isomorphism by $\mu:{\frak{S}}_g(\C) \to \CG$. 

For $A\in {\rm Mat}_g(\C)$, let $A^{ss}$ and $A^{nil}$ denote the semisimple 
and nilpotent parts, respectively. 
We note that $A^{ss}$ is also a symmetric matrix. 

\begin{lmm}
If $A\in \CH_g$, then the eigenvalues of $A^{ss}$ are all in $\CH$ and 
there is a complex orthogonal matrix $R$ such that $R^{-1}A^{ss}R$ is a diagonal matrix. 
\end{lmm}

\pr 
We define an inner product $(u,v)$ by ${}^tuv$ for $u,v\in \C^g$, and 
view $A$ as an endomorphism of $\C^g$.  
If an eigenvalue $\lambda$ of $A$ is not in $\CH$, then $A-\lambda I_g$ is again in $\CH_g$ and thus nonsingular. However, the determinant $\det (A-\lambda I_g)$ is zero, which is a contradiction. This proves the claim on the eigenvalues, and we now assume $A$ is semisimple. 
If $u,v\in \C^g$ are eigenvectors with different eigenvalues, say $\lambda$ and $\mu$, respectively, 
then since $\lambda {}^tuv={}^tu A v=\mu{}^t uv$, we have ${}^tu v=0$. 
Therefore, $\C^g$ is an orthogonal sum of eigenspaces as desired.
\prend
  
We also note that 
\begin{equation}
\CH_g^{ss}=\{ R^{-1}DR\in \CH_g \mid R\in O_g(\C), D \mbox{ is diagonal} \}  \notag 
\end{equation}
is a dense subset of $\CH_g$. Recall that every idempotent in $\CG$ is a conformal vector. 
Therefore, for each $A\in \CH_g^{ss}$ there are mutually orthogonal conformal vectors $e^1,\ldots,e^g$ with 
central charges $c_1,\ldots,c_g$, respectively, and 
scalars $\tau_1,\ldots,\tau_g$ such that 
$\sum_{j=1}^g e^j=\omega$ and $\mu(A)=\tau_1e^1+\cdots +\tau_ge^g$. 
We again note that $\tau_i\in \CH$.
Then $\{\widetilde{e}^j=(2\pi i)^2(e^j-\frac{c_{j}}{24}\1)\mid j=1,\ldots,g\}$ 
are mutually orthogonal conformal vectors for $(V,Y[\cdot],\1, \widetilde{\omega})$. 
Let $\wt^A_j[\cdot]$ denote the weight given by $\widetilde{e}^j$.  
Clearly, $\mu(-A^{-1})=\frac{-1}{\tau_1}e^1+\cdots+\frac{-1}{\tau_g}e^g$ for $A\in \CH_g^{ss}$.  
Therefore, by Theorem \ref{MMT}, we have the following result.

\begin{lmm}\label{SMT1} 
Let $A\in \CH_g^{ss}$ and $w\in V$ be a multi-$\wt^A_i[\cdot]$-homogeneous element. 
If $\widehat{{\rm Tr}}(w:A)$ is an analytic function on $\CH_g$, then for any $v\in \otimes_{j=1}^g{\rm VOA}(e^j)w$,  we have 
$$\widehat{{\rm Tr}}_{W^j}\left(o(v):-A^{-1}\right)
=\prod_{p=1}^g (\tau_p)^{\wt_p[v]}\sum_{h=1}^r s_{jh}\widehat{{\rm Tr}}_{W^h}\left(o(v):A\right).
$$ {\vspace{-3mm} \hfill \qed}
\end{lmm}

%Since $\wt^A_j[\1]$ is always zero, Theorem \ref{SMT2} is an immediate consequence of Lemma \ref{SMT1}.\\

Since $\wt_j^A[\1]$ is zero for all $A\in \CH_g^{ss}$, Theorem \ref{SMT2} is an immediate consequence of Lemma \ref{SMT1} for such $A$. Meanwhile, because $\CH_g^{ss}$ is a dense set of $\CH_g$
and $\widehat{{\rm \Tr}}_{W^j}(o(\1),A)$ and $\widehat{{\rm \Tr}}_{W^j}(o(\1),-A^{-1})$ are analytic for all $A\in \CH_g$ and $j=1, \dots ,r$, Theorem \ref{SMT2} holds as claimed.

%For $A\in \CH_g^{ss}$, since $\wt_j^A[\1]$ is alwasys zero, Theorem \ref{SMT2} for $A\in \CH_g^{ss}$
%is an immediate consequence of Lemma \ref{SMT1}.  Since $\CH_g^{ss}$ is a dense set of $\CH_g$
%and $\Tr_{W^j}(o(\1),A)$ and $\Tr_{W^j}(o(\1),Z^{-1})$ are all analytic on $\CH_g$, we have 
%the desired equation. 

%******************
%******************
%******************

\section{Applications}\label{Section-Applications}
As a corollary to Theorem \ref{SMT2}, we will recover the inversion transformation formula of matrices $A\in \CH_g$ for Siegel theta series. To do so we must introduce the lattice VOA $V_L$ for an even positive lattice $L$ of rank $g$.

We begin by first defining the VOA $M(1)$ of free boson type. Viewing $\C L$ as a $g$-dimensional vector space with inner product $\langle\cdot,\cdot\rangle$, 
we consider an affine Lie algebra 
\begin{equation}
\widehat{\C L}:=\left( \mathop{\oplus}_{\substack{j=1}}^g \mathop{\oplus}_{\substack{n\in \Z}}\C a_j(n)\right) \oplus \C, \notag
\end{equation}  
where $\{a_j\mid j=1,\ldots,g\}$ is an orthonormal basis of $\C L$ and 
the Lie products are given by $[a(n), b(m)]=\delta_{n+m,0}n\langle a,b\rangle$. 
We note $\widehat{\C L}$ does not depend on the choice of the orthonormal basis. 
Clearly, $\widehat{\C L}_{\geq 0}:=\left( \oplus_{j=1}^g\oplus_{n\geq 0}\C a_j(n)\right) 
\oplus \C$ is a subring. 
For every $\alpha\in \C L$, we define a one-dimensional $\widehat{\C L}_{\geq 0}$-module 
$\C e^{\alpha}$ by 
\begin{equation}
\begin{aligned}
a(n)e^{\alpha}&=0 \mbox{ for }n>0, \mbox{ and} \\
a(0)e^{\alpha}&=\langle a,\alpha\rangle e^{\alpha}. 
\end{aligned}
\label{App1}
\end{equation}
We also consider the induced module 
$$M^g(1)e^\alpha:=U(\widehat{\C L})\otimes_{U(\widehat{\C L}_{\geq 0})}\C e^{\alpha},$$
where $U(R)$ denotes the universal enveloping algebra of $R$.  
Among these modules, $M^g(1) e^0$ has a VOA structure of central charge $g$ which we denote 
by $M^g(1)$. Furthermore, 
$M^g(1)e^{\alpha}$ is an $M^g(1)$-module for each $\alpha$. Often $M^g(1)$ is 
called the VOA of $g$ bosons. 
Then 
$$V_L=\mathop{\oplus}_{\substack{\alpha\in L}} M^g(1) e^{\alpha}$$ 
becomes a VOA of central charge $g$, which is called a lattice VOA. 
(See \cite{FLM} for more details on lattice VOAs.) 
We note that $\1:=1\otimes e^0$ and 
$\omega:=\frac{1}{2}\sum_{i=1}^g a_i(-1)a_i(-1)\1$ are the Vacuum and Virasoro elements, respectively, of both $V_L$ and $M^g(1)$. 
  
It is known that $V_L$ is a regular VOA, and its simple modules are 
given by $V_{L+\beta}=\oplus_{\alpha\in L} M^g(1)e^{\alpha+\beta}$ for some 
$\beta\in \Q L$ with $\langle \beta, L\rangle\subseteq \Z$ (see \cite{D}). 
We will use the vertex operators  
\begin{equation}
\begin{aligned}
&Y(a(-1)\1,z)=\sum_{n\in \Z} a(n)z^{-n-1} \hspace{3mm} \mbox{  and}\cr
&Y(a(-1)b(-1)\1,z)=\sum_{m\in Z}\left(\sum_{n\in \N}a(-1-n)b(m+n)+b(-1+m-n)a(n)\right)z^{-m-1}.
\end{aligned}
\notag
\end{equation}
From (\ref{App1}), we have $o(a(-1)b(-1)\1)e^{\beta}=\langle a,\beta\rangle\langle b,\beta\rangle e^{\beta}$ and  
$\wt(a(-i_k)\cdots a(-i_1)e^{\alpha})=i_1+\cdots+i_k+\frac{\langle \alpha,\alpha\rangle}{2}$.  
Therefore, the character $\widehat{{\rm Tr}}_{M^g(1)}\left(o(\1):\tau\right)$ of $M^g(1)$ is $1/\eta(\tau)^g$ and 
the character of $V_L$ is $\theta_L(\tau)/\eta(\tau)^g$, where $\theta_L(\tau)$ is the theta series associated to the lattice $L$ and $\eta(\tau)=q^{\frac{1}{24}}\prod_{n=1}^{\infty}(1-q^n)$ is the eta-function, where $q=e^{2\pi i\tau}$.

By using an orthonormal basis $\{a_i\mid i=1,\ldots,g\}$ of $\R L$, 
we define $\omega^{ij}=\frac{1}{2}a_i(-1)a_j(-1)\1$ for $i,j=1,\ldots,g$. 
We note $\omega^{ij}=\omega^{ji}$, 
$\{\omega^{ii}\mid i=1,\ldots,g\}$ is a set of mutually orthogonal 
conformal vectors of central charge $1$, and 
$\omega=\sum_{i=1}^g\omega^{ii}$ is a Virasoro element of $M^g(1)$. 

From the construction, $M^g(1)$ has an automorphism $\sigma$ induced 
from $-1$ on $\C L$, that is, 
$\sigma(a_{j_k}(-i_k)\cdots a_{j_1}(-i_1)\1)=(-1)^ka_{j_k}(-i_k)\cdots a_{j_1}(-i_1)\1$.  
Let $M^g(1)^+$ denote the fixed point subVOA of $M^g(1)$ by $\sigma$. 
Then by direct calculations, we have $(M^g(1)^+)_0=\C \1^{\otimes g}$, $(M^g(1)^+)_1=0$, and 
$(M^g(1)^+)_2=\prod_{1\leq i\leq j\leq g} \C \omega^{ij}$ is isomorphic to a Jordan algebra of type $B_g$ 
by the $1$-products. 

We now introduce a multivariable trace function on the Siegel upper half-space $\CH_g$. 
For $A=(\tau_{ij})\in \CH_g$ and a $V_L$-module $M$, 
we recall the function (\ref{MultiTraceFct2}), and in particular 
$$\widehat{{\Tr}}_M\left(o(\1):A\right)={\rm Tr}_M e^{o\left(2\pi i \left(\mu(A)-\frac{\tr(A)}{24}\right) \right)},$$
where in this case, $\mu(A)=\sum_{i=1}^g\sum_{j=1}^g \tau_{ij}\omega^{ij}\in M^g(1)^+_2$. 

In order to pick out the lattice parts, we define 
$$\gamma_M(A)=\widehat{{\Tr}}_M\left(o(\1):A\right)\prod_{i=1}^{g}\eta(\mu_i)$$
for $A\in \CH_g$, where the $\mu_i$ are the numbers satisfying $\det(xE_g-A)=\prod_{i=1}^g (x-\mu_i)=0$. 
We note $\mu_i\in \CH$. \\
\indent We now prove the following result.

\begin{prn} \label{PropLattice}
For a lattice VOA $V_L$ with inequivalent simple $V_L$-modules $V_L=W^1, \dots ,W^r$, we have 
$\widehat{{\rm Tr}}_{W^h} \left(o(\1):A\right)$ is an analytic function on $\CH_g$ for all $h=1, \dots ,r$. Furthermore, 
$\gamma_{W^j}(A)$ are ordinary Siegel theta series and 
$$\left(-i\frac{1}{\det(A)}\right)^{g/2}\gamma_{W^j}(-A^{-1})=\sum_{h=1}^r s_{jh}\gamma_{W^h}(A),$$ 
where the $s_{jh}$ are the coefficients $A_{jh}^S$ in (\ref{ZhuThm}) for the matrix $S=\left(\begin{smallmatrix} 0&-1\\ 1&0\end{smallmatrix}\right)$. It follows that the functions $\gamma_{W^j}(A)$ satisfy the transformation laws of Siegel modular forms (see, for example, \cite{Freitag} for these equations).
\end{prn}

\pr 
As discussed above, a simple $V_L$-module $M$ is of the form $M=V_{L+\kappa}$ for some $\kappa\in \Q L$.
If $A=(\tau_{ij})$ is semisimple, then there is an orthogonal complex matrix $P\in O_g(\C)$ 
and scalars $\mu_1,\ldots,\mu_g$ such that 
$P^{-1}(\tau_{ij})P={\rm diag}(\mu_1,\ldots,\mu_g)$. 
Set $(b_1,\ldots,b_g)=(a_1,\ldots,a_g)P$ and $e^i=\frac{1}{2}b^i(-1)b^i(-1)\1$. 
Then $\{b_1,\ldots,b_g\}$ is an orthonormal basis of $\C L$, and $\{e^i\mid i=1,\ldots,g\}$ is 
a set of mutually orthogonal conformal vectors of $V_L$ such that 
$\mu(A)=\sum_{i=1}^{g} \mu_ie^i$. Since 
$$\widehat{{\Tr}}_M\left(o(\1):A\right)=\frac{1}{\prod_{i=1}^{g}\eta(\mu_i)}
\sum_{\beta\in L+\kappa}e^{\pi i\sum_{j=1}^g \mu_j\langle\beta,b^j\rangle^2}, $$
it follows that 
$$\gamma_M(A)=\sum_{\beta\in L+\kappa}e^{\pi i\sum_{j=1}^g \mu_j\langle\beta,b^j\rangle^2}.$$
Moreover, because $\pi i\sum_{j=1}^g \mu_j\langle \beta, b^j\rangle^2$ is an eigenvalue of 
$o(\pi i\sum_{j=1}^g \mu_j b^j(-1)b^j(-1)\1)$ 
for $e^{\beta}$, it is equal to an eigenvalue of 
$o(\sum_{j=1}^g\sum_{h=1}^g \tau_{jh}a^j(-1)a^h)$ for $e^{\beta}$, 
that is, 
$$\pi i\sum_{j=1}^g \mu_j\langle \beta, b^j\rangle^2 = \pi i\sum_{j=1}^g\sum_{h=1}^g\tau_{jh}\langle a^j,\beta\rangle\langle a^h,\beta\rangle.$$ 
Therefore $\gamma_M(A)$ is an ordinary Siegel theta series of $L+\kappa$. Explicitly, we have 
\begin{equation}
\gamma_M(A)=\sum_{\beta\in L+\kappa}
e^{\pi i\sum_{j=1}^g\sum_{h=1}^g\tau_{jh}\langle a^j,\beta\rangle\langle a^h,\beta\rangle}=
\sum_{\beta\in L+\kappa}e^{\pi i {}^t\widetilde{\beta}A\widetilde{\beta}}, \label{App2}
\end{equation}
where $\widetilde{\beta}=(\langle a^1,\beta\rangle,\ldots,\langle a^g,\beta\rangle)\in \R^g$. 
Since ${\rm Im}(A)$ is positive definite, there is a number $\epsilon (A) >0$ such that 
${}^t\widetilde{\beta}{\rm Im}(A)\widetilde{\beta}\geq \epsilon (A)
\langle \beta ,\beta \rangle$ for all $\beta\in L+\kappa$. 
It follows that 
$$\left|\gamma_M(A)\right|\leq \sum_{\beta\in L+\kappa}
\left| e^{-\pi \epsilon (A)\langle \beta,\beta \rangle}\right|<\infty.$$  
This implies $\gamma_M(A)$ is an analytic function 
for any symmetric matrix $A\in \CH_g$. 
Furthermore, (\ref{App2}) is well-defined for any $A\in \CH_g$, and so $\gamma_M(A)$ 
is an analytic function on $\CH_g$. \hfill \qed

\begin{rmk}
(i) Although we have been treating the cases where the rank $g$ of a lattice coincides 
with the genus of the Siegel upper half-space, 
by viewing $\CH_h\otimes I_{g/h}\subseteq \CH_{g}$ for $h|g$, we may treat a Siegel 
upper half-space of genus $h<g$. 

(ii) Let $\{W^1,\ldots,W^r\}$ be the set of simple inequivalent $V$-modules. Then
as Huang has proved in \cite{H}, 
$\sum_{i=1}^r \widehat{{\Tr}}_{W^i\otimes W^{i'}}\left(v:\tau\right)$ is invariant for 
an inversion $\tau \mapsto -\frac{1}{\tau}$, 
where $W^{i'}=\oplus_{p\in \C}{\rm Hom}(W^i_p,\C)$ denotes the restricted dual of $W^i$. 
Therefore, 
$\sum_{i=1}^r \widehat{{\Tr}}_{W^i\otimes W^{i'}}\left(o(\1):A\right)$ is invariant 
for an inversion $A\mapsto -A^{-1}$ by viewing $\CH_g\otimes I_2\subseteq \CH_{2g}$. 
\end{rmk}

\end{document}